\documentclass[11pt]{article}

\usepackage{amsmath, amsthm, amsfonts}

\usepackage[utf8]{inputenc}

\usepackage{graphicx}
\usepackage{xcolor}

\usepackage[T1]{fontenc}
\title{Mathematics \\ 
	an Imagined Tool for Rational Cognition}

\date{}

\author{Boris Čulina\footnote{University of Applied Sciences Velika Gorica, Department of Mathematics, e-mail: boris.culina@vvg.hr}}

\begin{document}
\maketitle

\noindent \textit{``...numbers are free creations of the human mind; \\
	they serve as a means of apprehending more easily and more sharply the difference of things.''}

\begin{flushright}
	Richard Dedekind \cite[p.~vii]{Dedekind}
\end{flushright}

%\textbf{Abstract.} Mathematics is the process and the result of shaping our intuitions and ideas   about our inner human world into thoughtful models  which are powerful tool in the process of rational cognition of the external world. Mathematical models are imagined worlds which can be only partially realized in our inner world. Consequently,  mathematical objects do not exist in reality and mathematical truths are not true in reality. However, in the process of rational cognition we synthesize mathematical objects and mathematical truths into truths about the real world.
 
\textbf{Abstract.} By analysing several characteristic mathematical models: natural and real numbers, Euclidean geometry, group theory, and set theory, I argue that a mathematical model in its final form is a junction of a set of axioms and an internal partial interpretation of the corresponding language. It follows from the analysis that (i) mathematical objects do not exist in the external world: they are imagined objects, some of which, at least approximately, exist in our internal world of activities or we can realize or represent them there; (ii) mathematical truths are not truths about the external world but specifications (formulations) of mathematical conceptions; (iii) mathematics is first and foremost our imagined tool by which, with certain assumptions about its applicability, we explore nature and synthesize our rational cognition of it.

\textbf{keywords:} mathematical intuition; mathematical models; mathematical objects; mathematical truths; applicability of mathematics

The basic problem of the philosophy of mathematics (not mathematics itself) is to answer the following intertwined questions:

\begin{itemize}
	\item Whether mathematical objects exist, and if so, in what way do they exist?
	\item What is the mathematical truth and how do we establish it?
	\item How is mathematics applied?
\end{itemize}

\noindent This article presents a solution that can be considered an elaboration of Dedekind's quotation cited at the beginning of the article. The basic theses that I intend to argue  are the following:

\begin{itemize}
	\item Mathematical objects do not exist in the external world. They are  imagined objects, some of which, at least approximately, exist in our internal world of activities or we can realize or represent them there.
	\item Mathematical truths are not truths about the external world but specifications (formulations) of mathematical conceptions. 
	\item 	 Mathematics is first and foremost our imagined tool by which, with certain assumptions about its applicability, we explore nature and synthesize our rational cognition of it.
\end{itemize}

\noindent I will try to make clear what is absolutely clear to the famous physicist  Percy W. Bridgman: ``It is the merest truism, evident at once to unsophisticated observation, that mathematics is a human invention.'' \cite[p.~60]{Br}.  Having practised mathematics all my life, by vocation and by profession, just as breathing was natural to me, so it was natural for me to consider mathematics as a human invention and a free creation of the human mind whose purpose is to be a tool for our rational cognition and rational activities in general.\footnote{Of course, many mathematicians do not share my opinion that mathematics is a free creation of the human mind.} When I decided to clarify to myself what  human invention mathematics was, it proved to me, I believe, a far more difficult task than explaining what breathing is. In this article, I have outlined what I came up with along the way. 

In most of the article I will explain and argue the view of mathematics as an imagined tool for rational cognition. In sections 1 to 6 I will consider as illustrative examples the classical mathematical models which are still the most important ones: the natural number system, the real number system, and Euclidean geometry; as well as today’s standard mathematical models: group theory, and set theory (sections 1 to 6). To quote the famous mathematician Saunders Mac Lane: ``\ldots a philosophy of Mathematics is not convincing unless it is founded on an examination of Mathematics itself.'' \cite[p.~60]{Lane}. I will then set out the basic characteristics of the view of mathematics as an imagined tool for rational cognition (Section 7), explain how mathematical objects possibly exist (Section 8) and how mathematics is applied (Section 9). Section 10 describes how the whole structure of mathematics can be understood as an imagined tool for rational cognition. In the last part of the article I will expose  this  view of mathematics to Benacceraf dilemma \cite[p.~661]{Be}, Bueno's five desiderata \cite[p.~63]{OB} and  Quine-Putnam indispensability argument \cite{Co}  (Section 11), as well as  determine its place on the map of different views of mathematics (Section 12).

The terms ``intuition'', ``idea'', ''conception'', ''model'' and ''theory''  will denote more and more precise stages in modelling  certain thoughts. The terms ``model'' and ``theory'' will eventually acquire a precise meaning that deviates from the standard meaning in mathematical logic. They will be synonymous here and will denote a set of axioms in a language, together with a  partial interpretation of the language. In doing so, I will give preference to the term ``model'' when the emphasis is on the interpretation and to the term ``theory'' when the emphasis is on axioms.  Thereby, the term ``interpretation'' generally refers to the meanings of linguistic forms, and more specifically to the interpretations of first-order languages. The term ``specification'' generally refers to the refinement of ideas, and more specifically to the description using a set of sentences (axioms) of a first-order language that an interpretation of the language must satisfy. By ``internal world of activities'' I mean the world that  consists of activities over which we have a strong control, and which we organize and design by our human measure. For example, these include movements in a safe space, grouping, arranging and connecting small objects, spatial constructions and deconstructions with small objects, talking, writing and drawing on paper, shaping and transforming manipulative material, combining and repeating actions, making choices, dynamics of actions and changes in the environment subordinated to us, painting, singing, etc. This does not include activities with objects over which we have no strong control, or our activities are significantly limited by the environment, such as e.g. climbing steep rocks, building a house, etc. Although each person has her own specific way of organizing and designing internal activities, which depends on her as well as on the environment to which she belongs, due to the generality of further considerations, it does not matter whether we mean the internal world of activities of an individual person or the internal world of activities of the entire human species. The internal world of activities will be discussed in more detail in Section 7. I would also note here that all imagined objects are considered as concrete and not abstract objects. In Section 7 it will be argued that language is the bearer of abstraction and not objects.

%By ``external world'' I will mean the world that would remain even if we became extinct as a species (e.g., books would remain but not the meaning of what is in the books), but also the world that would disappear with us and that is not under our direct control, (e.g., our biological tendency to dance but not a concrete performance).

% some standard topics and two seemingly close views, structuralism and fictionalism, in modern philosophy of mathematics.

 \section{Natural numbers}

Natural numbers are the result of  modelling our intuition about the size of a collection of objects. We measure the collection through the process of counting, and natural numbers are objects for counting. To start counting we must have the first number, to associate it with the first chosen object in the collection. To continue counting, after each number we must have the next new number in order to associate it with the next chosen object in the collection. There is no special reason to sort out certain particular objects as natural numbers. Merely for the needs of performing a calculation  we sort out a particular realization.  In the past those were collections of marbles on an abacus, and today we use sequences of decimal digits on paper and of bits in a computer. It means that for counting it is not important how numbers are realized, but only the structure of the set of natural numbers which enables us to count is important: that there is a first number and that each number has its successor.  It seems that natural numbers exist in the same way as chess figures, in the sense that we can always realize them in some way. However, the structure of natural numbers, as opposed to the structure of chess, brings in itself an idealization. So that the counting could always continue, each natural number must be followed by the next natural number. Therefore, there are infinitely many natural numbers. Thus, although we can say for small natural numbers that they exist in the standard sense of that word, the existence of big natural numbers is at best a kind of  idealized potential existence.  

In order to precisely formulate the conception of natural numbers, we need a corresponding sufficiently precise language, a mathematical language, or even better, a first-order language.  Among other symbols, the language should contain the symbol  for the first number, the standard is ``1'', and the function symbol, for example an ``S'', for the immediate successor operation ($n\mapsto \textrm{S}(n)$). This allows us to name each number. What  the names denote  is not so important. It is only essential that the named objects have the structural role of natural numbers (the names can even denote themselves).  We specify the conception by declaring certain sentences of the language to be true about natural numbers. It is not necessary to precisely specify these claims here, nor the language in which they were made.\footnote{The standard formulation consists of Peano's axioms and recursive conditions for addition and multiplication in the first-order language of  arithmetic \cite[p.~82]{vanD}.}   In that language, the claims must express, inter alia, (i)  that 1 is the first number, (ii) that each number $ n $ has its successor $ \textrm{S}(n) $ which is a new natural number in relation to all previous natural numbers, and (iii) that any natural number can thus be obtained. I will hereinafter call these claims the axioms of natural numbers. In my view, the axioms of natural numbers are neither true nor false, just as the axioms that would describe the game of chess would be neither true nor false. They are simply a means of specifying our ideas about the objects we use for counting.\footnote{In \cite[Chapter~7]{FerP} Ferreirós argues that Peano's axioms are self-evident truths about our practice of counting. Our practice of counting is not a natural process but a voluntary activity. That is why, in my opinion, the truth of Peano's axioms does not derive from that practice, but from the fact that they norm that practice, especially since they introduce an element of idealization into the practice.} However, G\"{o}del's  incompleteness theorems \cite{G} and Lowenheim-Skolem theorems \cite{L,S} tell us that we cannot have a complete (at least up to isomorphism) specification in a first-order language, even if we include sets in the specification.\footnote{Although in the classical set theory all structures of natural numbers are mutually isomorphic, according to the Lowenheim-Skolem theorems, the classical set theory itself has non-isomorphic interpretations \cite[p.~105-106]{vanD}.} Therefore, in addition to the axiomatic specification, it is necessary to have an interpretation of the language, as well. In the case of natural numbers, it is a partial internal interpretation --  a partial interpretation in our internal world of activities. The interpretation is partial due to the idealization of the existence of extremely big numbers. The interpretation belongs to our internal world of activities because numbers are our imagined constructions that we can partially, to the point of isomorphism, realize in the world available to us: using marbles on an abacus, tallies on paper, etc.  On the other hand, even in the idealization of the total interpretation of the language, it is necessary to additionally have an axiomatic specification, because the recursively defined truth value of sentences in the interpretation is not a computable function \cite{Tarski, G}. Since we can never construct all numbers, the overall structure of natural numbers does not exist in the literal sense of the word. There is only the conception of natural numbers, specified by axioms, and  partially realized in our internal world of activities.  This is the final result -- mathematical model -- of modelling our intuition about natural numbers as objects for counting. It carries with it incompleteness and the ever-present tension in mathematics between basic intuition and the constructed model, as well as between the axiomatic specification and the constructive content of mathematical concepts. This tension is a positive source of accepting new axioms as well as improving intuition.\footnote{This tension took a more dramatic form in the historical conflict between modern and classical mathematics, primarily in the confrontation between Dedekind's modern conceptualist approach and Kronecker's classical constructivist approach, as well as in the crisis in the foundations of mathematics -- in the conflict between Hilbert's program in the foundations of mathematics and Brouwer's intuitionism (see, for example, \cite{FerC}). How much of this tension is present in the Hilbert's program itself can be read in \cite{Sieg}.}

\section{Real numbers}
 Through real numbers we organize and make precise our intuition about the process of measuring.\footnote{In my opinion, this is the main role of real numbers, as a tool for rational cognition. Various practices with real numbers and their central role in mathematics are nicely described in \cite[Chapter 8]{FerP}.}  While natural numbers are imagined as objects for counting, real numbers are imagined as the results of the process of measuring.  However, we can imagine (idealized) situations in which the process of measuring never stops -- we generate a potentially infinite list of digits, with no consecutive repetition of the same group of digits after an individual step. If we want to have the results of such processes of measuring, we must introduce, in addition to rational numbers,  new results of measuring -- irrational numbers.  As opposed to natural numbers whose existence we can understand at least as some kind of an idealized potential existence, we cannot explain the existence of irrational numbers in this way. Although we can approximate irrational numbers by rational numbers with arbitrary precision, their existence is outside  our means of construction -- we have just imagined irrational numbers.

 As is the case with natural numbers, the final mathematical model of real numbers is a junction of axiomatic specification\footnote{The axiom of completeness ensures that any idealized measurement process has a result. A variant of the axiom corresponding to this approach postulates that every decimal expansion $ a_0, a_0.a_1, a_0.a_1a_2, \ldots $ (where $ a_1,a_2,\ldots $ are decimal digits, whereas $ a_0 $ is an integer) has a limit and that every real number is a limit of such an expansion. All other variants of the completeness axiom postulate in other ways these imagined objects. See, for example, \cite{Teis} -- the variant with decimal expansions has the label CA18.}  and partial internal interpretation of the corresponding language. For example, in the mathematical model we can identify Euler number  $e$, the irrational number to whom  the sequence  $\Big (1+\dfrac{1}{n}\Big )^n$ is  closer, as we increase the natural number $n$. Although we can approximate the number $e$ with arbitrary precision by constructions in our internal world of activities,  it certainly does not exist in the same way as  my dog exists. It exists in the same way as an idealized material point in classical mechanics, as a non-existing phlogiston in a wrong theory about chemical reactions, and as Snow White in the classical fairy  tale \textit{Snow White and the Seven Dwarfs}. However, although our language usually has only a partial interpretation, the classical logic of using the language assumes that it is a semantically complete language -- that it has a complete interpretation: each name names an object, each predicate symbol refers to a binary predicate, and each function symbol refers to a function.\footnote{For complete interpretations of the first-order languages see, for example, \cite[Section~3.4]{vanD}}  Because of this assumption, in thinking itself there is no difference whether we think of objects that really exist or we think of objects that do not really exist. That difference can be registered only in a ``meeting'' with reality. And for mathematics there is no such meeting: a mathematical model creates its own reality in our internal world of imagination. However, unlike erroneous physical models, the mathematical model of real numbers can be realized approximately in  our internal world of activities in the same way that correct physical models are realized approximately in the external world or children’s fairy tales in real theatrical performances. 

\section{Euclidean geometry}

With the appearance of non-Euclidean geometries in the 19th century, Euclidean geometry lost the status of an a priori mathematical theory. It became only one of the possible models for physical space, distinguished only by the fact that it is a good approximation of the space in which practical science takes place.\footnote{In \cite{T}, the impact of the appearance of non-Euclidean geometries on the philosophy of geometry is analyzed in detail.} Contrary to such a view, according to which Euclidean geometry is a part of physics, I will argue here that just as number systems are idealized conceptions derived from intuition about our internal activities of counting and measuring, so too is Euclidean geometry an idealized conception derived from intuition about our internal spatial activities. Since the view of Euclidean geometry presented here is not standard, the argumentation will be more detailed than in the sections above.

Internal spatial activities should be distinguished from external spatial activities. The former are conditioned by our human nature, the latter additionally by the world around us. At first glance, it seems difficult, almost impossible, to draw a clear line between these two types of activities. However, for some activities we can clearly determine that they belong to external spatial activities. For example, mountain climbing is an external spatial activity because it involves orientation in a given landscape and taking care of the configuration of the terrain over which we move.  When we are on the different parts of a mountain road, different spatial situations will require different responses. However, as far as our ability to react is concerned, it is the same in all places and in all directions. Of course, there is also the ubiquitous gravitational force that makes one direction in space prominent, which we especially have to take care of. Likewise, if there is a strong wind blowing from a certain direction, that direction also becomes prominent to us. However, we always attribute the appearance of differences in a certain direction to the influence of an external factor, which shows that a priori all directions are the same to us. If we have to light a fire by placing twigs so that they form a cone,  our approach to geometric construction will be the same, whether we are making a small or a large cone. This shows that our ideas of spatial constructions are independent of  the units we use in their construction. If for a certain selection of units a change in  construction occurs, we ascribe it to an external factor. If in this way we try to identify the nature of our internal spatial activities, as invariants to the different spatial situations in which we find ourselves, then we are very close to Delboeuf's analysis \cite{Delboeuf}. He considers what remains when  we ignore all differences of things caused by their movements and mutual interactions. According to Delboeuf, in the ultimate abstraction from all diversities of real things we gain the homogeneous (all places are the same), isotropic (all directions are the same), and scale invariant (geometric constructions are independent of size) space -- the true geometric space which is different from the real space. However, for Delboeuf this geometry is the background geometry of real space, while for me it is the geometry of our internal activities in space. 

We can come to the same conclusion if instead of an external argument, seeking common  ground in all our external spatial activities, we use an internal argument, analysing our internal spatial activities directly and independently of the external world. A simple introspection  shows that we do not distinguish different places, different directions and different units for spatial constructions until the outside world forces us to distinguish them. To eliminate the presence of gravity on the Earth's surface we must look for examples where it is negligible. In addition to the extravagant situation of a free fall, these can be examples of activities that take place approximately in the horizontal plane or three-dimensional examples in which gravity is not important.  For example, a child will build his imaginary monster using Lego bricks regardless of the location where he built it, the way he oriented it in space and the dimensions of basic Lego bricks used in its construction. The same indifference to location, direction, and size is present when we rearrange the Rubik's cube. 

I believe that our most basic approach to space, the approach inherent to us, is an  \textit{a priori ignorant  approach to space}: all places are the same to us (the homogeneity of space), all directions are the same to us (the isotropy of space) and all units of length we use for constructions in space are the same to us (the scale invariance of space). These three principles of symmetry express our basic intuition about our internal spatial activities. Any deviation from these symmetries we attribute to the external world. Thus, it is precisely these principles of symmetry that determine a clear boundary between our internal and external spatial activities.\footnote{The importance of these principles has been recognized a long time ago in the works of Delboeuf \cite{Delboeuf}, Helmholtz \cite{Helmholtz}, Clifford \cite{Clifford1, Clifford2} and Poincar\'{e} \cite{Poincare}, but in a different interpretation than the one described here.  In the 17th century John Wallis  proved, assuming other Euclid's postulates, that  the scale invariance principle saying that ``For every figure there exists a similar figure of arbitrary magnitude.''  is equivalent to the Euclid's fifth postulate \cite{Wallis}. Also, it is well known that among all Riemannian manifolds Euclidean geometry is characterized by these three symmetry principles (see, for example, \cite{Clifford1}).}

Just as internal and external spatial activities can be distinguished, so can the geometry resulting from these activities. In my view of mathematics as our, in a certain sense, a priori tool of rational cognition,\footnote{I will explain my views below in Section 7.} geometry arising from internal activities is a part of mathematics, while geometry arising from external spatial activities is a part of physics. For example, caring about the direction of gravity belongs to external spatial activities, so pointing out that direction belongs to physical geometry and not to geometry arising from our internal spatial activities. Of course, in this physical geometry, as in other physical theories, appropriate mathematics is incorporated, but it is a different tool of rational cognition (see Section 5) than the geometry that arises from our internal spatial activities. In the book \cite{Sche}, it is nicely described how appropriate physical geometries emerge from our external spatial activities, which are very important for our survival.

In \cite{Cu3}, an elementary system of the axioms of Euclidean geometry is developed. On the one hand, the system is directly founded on the three principles of symmetry described above, while on the other hand, through the process of algebraic simplification, it gives an equivalent Weyl's system of axioms of Euclidean geometry (the axioms of Euclidean affine space) \cite[Chapter 1]{Weyl}. In this way, Euclidean geometry is characterized by these three principles of symmetry without any additional assumptions (except the idea of continuity). Thus, I gave the argument that Euclidean geometry is an idealized mathematical model derived from intuition about our internal spatial activities.

I consider this interpretation in space of our human internal activities the primary interpretation of Euclidean geometry.  However, we can preserve the sentence part of the theory but change the interpretation. Then it does not need to be a mathematical conception anymore. It depends  on a new interpretation, be it an external  or an  internal one. If we ask ourselves whether the physical space obeys the axioms of Euclidean geometry, we must extract from  space what  we consider as points (maybe enough localized parts of space), as directions (maybe  directions of light rays), and   the distance between two points (maybe the time needed for light to pass from one point to another). If in  such an interpretation the physical space satisfies the axioms of Euclidean geometry then we have an experimentally verifiable theory. Its sentence part is the same as in our mathematical theory of the space of our human activities, so we can transfer all results to the structure of  physical space. Only the interpreted part is different. It does not belong to mathematics anymore, but it is a base for an experimental verification of the theory about the external world. However, we can change an interpreted part of the originally imagined Euclidean geometry in a way that it will be still a mathematical theory. And it happens in mathematics  often. Namely, when we investigate complex mathematical objects which  we cannot perceive so easily, for example a set of functions of some kind, it is useful to find Euclidean structure in it. Then we can transfer our geometric intuition to that set -- think of functions as points,  measure how distant two functions are, etc. In that way we can visualize them and succeed in thinking about them more easily and  effectively.

The example of Euclidean geometry purports  that only the axiomatic part of a theory can belong to mathematics, while the interpretation does not have to. Also, mathematical interpretation does not have to be an idealized direct interpretation in our internal world of activities, but it can also be an interpretation in another mathematical model (theory).

\section{Group theory}
Group theory, in addition to being an elementary part of more complex mathematical theories, models above all our intuition about symmetry as invariance to certain transformations. Since different situations have different symmetries, unlike previous mathematical models which have an intended interpretation, this theory does not have an intended interpretation, but has intended non isomorphic interpretations. Thus, some mathematical models are simply sets of axioms without a specific interpretation. However, if they are modelling some important inner intuition about our approach to the world, as group theory does, then they are usually very important. Probably, the most famous example is the Riemann's conception of geometry as a manifold with a metric \cite{R}. These models found their application half a century after their invention with the appearance of Einstein's general theory of relativity. Today, manifolds are an essential component of mathematics and physics. Although the application was realized so late, it had to happen, because manifolds model successfully the basic mathematical idea about the coordinatization of investigated objects, an idea that generalize such an efficient idea of measuring. Although, due to their generality, the theory of groups and the theory of Riemannian manifolds have no interpretation in our world of internal activities, they grew out of intuition about the world of internal activities, the former on the ideas of transforming objects and combining transformations, the latter on the idea of coordinatization of objects.

In the language of group theory, due to the existence of  non-isomorphic interpretations, we sometimes think of a definite, and sometimes of an indefinite (``any'') interpretation. However, the very use of language and its logic requires that when we think in a language, we necessarily assume that it is a semantically complete language, no matter how we imagine the interpretation. The situation is the same as when we use  variables in our thinking. Whether we attach a certain value to a variable or not, in thinking within a classical mathematical language we necessarily assume that it has a certain value. However, if we think of groups in the language of set theory, then the groups themselves are the values of variables, not interpretations of the whole language as described above, and we think of them differently. The language of set theory allows us to connect and compare groups with each other, without having to know the true nature of individual groups, but possibly their isomorphic copies in the world of sets without urelements. Thus, the process of modelling the initial intuition and the way of working with the constructed models depends on the language in which we model the intuition.

\section{Mathematical models from other disciplines}
  The source of mathematical models does not have to be an intuition about our internal world of activities. They can also be ``borrowed'' from other disciplines. The nature of our thought and use of language, as well as the way in which we manage  the vast complexity of the world, leads to the extraction of  a certain structure from such a domain. We extract from the domain certain objects, relations and operations and we describe their properties.   If we have thus obtained an important model from that domain then its sentence (axiomatic) part is a mathematical model important for examination. We can use classical mechanics to illustrate this point. Although particles, motion and forces do not belong to mathematics, mathematics can take the structural properties of phenomena (usually described as a set of sentences in an appropriate language) and formally investigate them: the consequences (for example, in the problem of three bodies), the equivalent formulations (for example, Lagrangian and Hamiltonian formulations of Newtonian mechanics emerged in this way), etc.
 
 \section{Set theory}
In the consideration of any objects, the consideration of the sets (collection) of those objects naturally occurs. In mathematics, this step has a deeper meaning. Namely, the foundational mathematical modelling must model  the very intuition about mathematical and thought modelling itself. In the process of thought modelling, we extract a structure from a set of objects, that is to say, we extract some distinguished objects, and some relations and functions over the set of objects. Therefore, the subject of the foundational modelling must be the structures themselves and their parts. We can reduce the description of the structures to the description of their parts. It is the standard result of mathematics that we can describe distinguished objects by functions, functions by relations, and relations by sets.  In this way we can reduce the foundational modelling to the analysis of sets. From sets we can build all structures. Also, we can compare such structures using the functions between them. The language of sets provides simple means for describing structures and constructing new structures from the old ones. In this way, sets give us an universal language for mathematics. Furthermore, sets are often necessary for specification. For example, the specification of natural numbers requires the axiom of induction, which, in its most general formulation, needs the notion of a set. Likewise, the specification of real numbers requires the axiom of continuity, and the specification of Euclidean geometry requires Hilbert's axiom of maximality, and they both need the notion of a set. However, in what way exist the set of natural numbers, the set of real numbers, and the set of space points, when they are infinite? Moreover, when we think of sets, we also consider sets of sets. If we want to have an elegant, rounded and universal set theory, infinite sets are naturally imposed on us, truly the whole infinite hierarchy of infinite sets, together with infeasible operations on them.\footnote{An elegant exposition of the standard ZFC set theory can be seen, for example, in \cite{End}.}  How can we understand the existence of such sets and operations? Should we reject this theory, which has proven to be very successful, because its objects can be realized only when they are finite? Hilbert, who certainly knew what good mathematics is, said the following on the Cantor's set theory of infinite sets:\footnote{Cantor's exposition of his set theory can be seen in \cite{C}.} ``This appears to me to be the most admirable flower of the mathematical intellect and in general one of the highest achievements of purely rational human activity,'' \cite[p.~167]{Hi}.\footnote{The complicated and at times dramatic historical development of set theory is described in detail in \cite{FerS}.}

The language of set theory presupposes an intended interpretation. In addition to the fact that we can only partially realize it, the interpretation itself is not clear to us in many ways. It is clear that the idea of a set derives from our activities of grouping, arranging and connecting objects and that set theory is an idealized mathematical model for these activities. However, in the very finite part, when we talk, for example, about the set containing three concrete objects $ o_1 $, $ o_2 $ and $ o_3 $, it is not clear what kind of object the set itself is, let us call it $ s = \{ o_1, o_2 , o_3 \} $. Formally, we can describe the situation by saying that we have added a new object to the objects, which we then call the set of these objects and which has the unique property that only the objects $ o_1 $, $ o_2 $ and $ o_3 $ and no other belong to it. Such a description would correspond to a combinatorial approach and obviously has a structuralist overtone -- the very nature of  sets is not important but their relationship to other objects is important. However, another description $ s = \{ x | x = o_1 \textrm{ or } x = o_2 \textrm{ or } x = o_3 \} $ of the same set has a different connotation. Now the set is given by a predicate, so it is a kind of extensional abstraction of an one-place predicate. In this view of sets, as extensions of one-place predicates, they have no structural role but have their own individual nature in our world of meaningful linguistic forms, in the same way that points and directions have their own individual nature in  our world of internal spatial activities. Let us note that both the structural and individual view of sets change the initial intuition about the impossibility of realizing infinite sets. In the structural view, the set of all natural numbers is simply a new object. It differs from a finite set only in that there are infinitely many objects in the membership relation with it.  Natural numbers can be defined without the concept of infinity.\footnote{see, for example, \cite[p.~75]{Q}} Thus, from a structural point of view, the formation of the set of all natural numbers is not burdened by the concept of infinity, although it has infinite members.  In the individual view, the predicate ``to be a natural number'', as commented above, can be defined without the concept of infinity. Let's compare, for example, that predicate with the predicate ``to be an electron''. The letter predicate can be defined by a certain experimental procedure. It is well defined regardless of whether there are finitely or infinitely many electrons in the world. It is the same with the predicate ``to be a natural number'': it is well defined regardless of whether there are finitely or infinitely many natural numbers. In the individual point of view, that predicate, as well as any other that has the same extension, can be considered the set of natural numbers. So, we can conclude that the formation of the set of natural numbers is not burdened by the concept of infinity, although it has infinite members. Thus, both in the structural and in the individual view, the set of natural numbers  is a finitely well-formed object. The idea that a set must be ``made'' of its elements is not present here at all, an idea according to which we can never make an infinite set. Despite all the doubts related to the notion of a set, the constructed mathematical model  is very successful. Today, ZFC axioms form its axiomatic part and the model has an intended  interpretation, although there are doubts as to what the interpretation is, what we can realize  and how we can realize it. Here we only have a more pronounced tension between the basic intuition and the final model, which can ultimately lead to model refinement, model change, or even separation into multiple models.\footnote{In addition to other models of set theory that complement, weaken or change the ZFC axioms (see \cite{Fra} and, for newer alternatives, \cite{Hin}), the multiverse view in set theory has recently been developed, according to which there is no ``true'' mathematical model for the concept of a set, but there are various equally acceptable models \cite{Ham}.}

 \section{What mathematics is}
The conceptions described above possess all of  the essential characteristics of mathematical conceptions. First of all, mathematical conceptions have a clearly defined purpose -- to be  successful tools  in the process of rational cognition, and in rational activities in general. This purpose significantly influences their design and determines their value. We can use mathematical conceptions  directly, like the use of numbers, through an ideal model of interaction with the world. We can use them indirectly: (i) like the use of the Euclidean geometry -- by changing the interpreted part of the theory into an external interpretation, (ii) like the use of the group theory -- by giving just the axioms and their consequences regardless of interpretation, which could be the external one, or (iii) like the use of the set theory -- by organizing effectively other mathematical tools. Also, we use mathematical conceptions indirectly, (iv) as constituent parts of more complex mathematical conceptions -- as it is, for example, the case with Euclidean space as a tangent space on a Riemannian manifold, or (v) we use them indirectly in the way described above with Euclidean geometry -- to interpret them in collections of complex mathematical objects for the purpose of making them more intuitive and more manageable. A multitude of specific mathematical models that are used to model specific problems should also be mentioned here. Such a model is applied directly, its purpose is concrete, and its design and evaluation largely depend on the problem it models.\footnote{In the book \cite{Len}, mathematics is analysed as a means of building specific mathematical models, while in the book \cite{Fen} a general model of such use of mathematics is given.}

%like  other useful types of conceptions that do not pretend to be  about the external world, have a goal. Fairy tales   have the goal to edify. The game of chess has  the goal  of being an intellectual amusement. The goal of mathematical conceptions is to be  successful tools  in the process of rational cognition, and in rational activities in general. And that goal determines their nature. We can use mathematical conceptions  directly like the use of numbers, through an ideal model of interaction with the world. We can use them indirectly: (i) like the use of Euclidean geometry -- by changing the interpreted part of the theory into an external interpretation, (ii) like the use of group theory -- by giving just the axioms and their consequences regardless of interpretation, which could be the external one, or (iii) like the use of set theory -- by organizing effectively other mathematical tools. Also, we use mathematical conceptions indirectly, (iv) as ingredients of more complex mathematical conceptions -- as it is, for example, the case with Euclidean space as a tangent space on a Riemann manifold, or (v) we use them indirectly in the way described above with Euclidean geometry -- to interpret them in collections of complex mathematical objects for the purpose of making them more intuitive and more manageable.

The conceptions described above purport that mathematics is an inner organization of rational cognition and knowledge, a thoughtful shaping of the part of the cognition that belongs to us.  For example, we organize the possible results of measurement into an appropriate  number system. The inner organization  needs to be  distinguished from (but not opposed to) the outer organization of rational cognition, a real shaping of an environment that  comprises construction of a physical means for cognition (for example, an instrument for  measuring  temperature). Mathematics is a process and a result of shaping our intuitions and ideas   about the reality of our internal activities, into thoughtful models which enable us to  understand and better control  the whole  reality. For example, we shape our sense for quantity into  a system of measuring quantities by numbers. Thoughtful modelling of other intuitions about our internal human world of activities, for example intuitions about symmetry, flatness, closeness, comparison, etc., leads to other mathematical models. This claim will gain its full meaning only when I explain what the terms ``internal world of activities'', ``intuition'' and ``mathematical model'' mean. That is the content of the rest of this section.

In the introductory part, I briefly described our internal world of activities as the world that  consists of activities over which we have a strong control, and which we organize and design by our human measure. The examples I listed there are activities that we can easily recognize in free children's play.\footnote{In \cite{CuW}, it is argued that, contrary to the narrow standards of mathematics education, we best help children in their mathematical development by providing them with an environment in which they will, in free play, and with our unobtrusive help, develop their internal world of activities, design it, conceptualize it and apply it to problem-solving.} These activities form the basis of the adult world of internal activities. These activities develop as we grow up, but this is primarily the development of their design and conceptualization. Their presence in early development stages indicates their biological basis.\footnote{The book \cite[p.~28]{Lakoff} discusses: ``ordinary cognitive mechanisms as those used for the following ordinary ideas: basic spatial relations, groupings, small quantities, motion, distributions of things in space, changes, bodily orientations, basic manipulations of objects (e.g., rotating and stretching), iterated actions, and so on.''} They are a unique characteristic of the human species, an essential part of our evolutionarily developed abilities by means of which, unlike other species that adapt to the environment, we adapt the environment to ourselves.\footnote{``Man is a singular creature. He has a set of gifts which make him unique among the animals: so that, unlike them, he is not a figure in the landscape –- he is a shaper of the landscape.'' \cite[p.~19]{Bro}.} Of course, cultural evolution and social context play a decisive role in designing and conceptualizing these activities.\footnote{See, for example, \cite{Kit}  for the influence of cultural evolution and  \cite{Hersh,Ernest} for the influence of human society.} 

There are many examples of highlighting the importance of these activities. Dedekind talks about ``the ability of the mind to relate things to things, to let a thing correspond to a thing, or to represent a thing by a thing, an ability without which no thinking is possible'' \cite[p.~viii]{Dedekind}. Hilbert sees the source of his finitist mathematics in ``extralogical concrete objects that are intuitively present as immediate experience prior to all thought'', and these objects are ``the concrete signs themselves, whose shape \ldots is immediately clear and recognizable'' \cite[p.~171]{Hi}. Feferman writes that the source of mathematical conceptions ``lies in everyday experience in manifold ways, in the processes of counting, ordering, matching, combining, separating, and locating in space and time'' \cite[p.~75]{Fef}. Hersh writes: ``To have the idea of counting, one needs the experience of handling coins or blocks or pebbles. To have the idea of an angle, one needs the experience of drawing straight lines that cross, on paper or in a sandbox'' \cite[p.~46]{He}. And so on. 

Internal activities are concrete activities with concrete objects: they take place in space and time, and in a given environment -- on the table with a pencil and paper, in a ballroom, or on a sandy beach. They are experiential activities, but it is not an experience of the external world, but an experience of our actions in the world available to us and subordinated to us. We experiment not with the objects of the external world but with the possibilities of our actions in the world adapted to us. In this world, Piaget distinguishes two types of knowledge \cite[p.~16-17]{Pia}. An example of the first type of knowledge is when we pick up two objects and determine that one of them is heavier. This knowledge arose from our action performed on the objects and its source is in these objects: it is knowledge about these objects -- knowledge about the world. Piaget calls such knowledge ``physical knowledge''. But when we order objects, for example by weight, this order was created by our actions and its source is in our activities and not in the objects. The knowledge that has its source in such activities, for example that there is always the same number of objects no matter how they are ordered, is called by Piaget ``logical mathematical knowledge''. Mathematical knowledge has its source in the world of internal activities precisely in this way: it springs from these activities themselves and not from the objects of these activities. Feferman writes: ``Theoretical mathematics has its source in the recognition that these processes are 
independent of the materials or objects to which they are applied and that they are 
potentially endlessly repeatable.'' \cite[p.~75]{Fef}. However, we cannot completely separate activities from the objects of the activities: generally speaking, they depend on the objects. For internal activities, it is important that this connection is weak. Our control over these activities is the dominant force, rather than the influence exerted over them by the objects.  

Likewise, it is not possible to draw a clear line where the world of internal activities ends, and activities become external. An example can be made out of constructing and deconstructing objects. When a child does this with Lego blocks, it is surely an internal activity. Although the structure of Lego blocks affects the possibilities of construction, these are activities over which we have a strong control and the possibility of designing them according to our own measure. Constructing and deconstructing stone walls without mortar is certainly not an internal activity, because it requires experience working with weights, centres of gravity and forces in contact. But both activities are the source of the same mathematical idea, the idea of analysis and synthesis of what we are researching: let us examine the phenomenon by breaking it into parts, study those parts and synthesize the knowledge thus acquired into knowledge about that phenomenon. This idea is a mathematical idea because it has its source in our approach to the world, not in the world itself. Although it is present in both mentioned activities, in the internal activity it takes a clear and separate form, while in the other activity it is connected with the physical content. Because of the freedom we have in the world subordinated to us, I believe that we can always internalize external activities, represent them with internal activities in which the mathematical idea will come to its full expression. 

To conclude, although it is about concrete activities in a real environment, due to the subordination of that environment to our activities and due to the strong control over these activities and the strong possibilities of their shaping, we can talk about these activities as our internal world of activities, and even about a certain kind of a priority of that world.

The concept of intuition is not an intuitive concept. In \cite[p.~49]{Kit} Kitcher writes: ``'Intuition' is one of the most overworked terms in the philosophy of mathematics, Frege's caustic remark frequently goes unheeded: "We are all too ready to invoke inner intuition, whenever we cannot produce any other ground of knowledge." ''.\footnote{For example, Brouwer in his First act of intuitionism \cite[p.~509]{Brou}writes: ``\ldots intuitionist mathematics is an essentially languageless activity of the mind having its origin in the perception of \textit{move in time}, i.e. of the falling apart of a life moment into two distinct things, one of which gives way to the other, but is retained by memory. If the two-ity thus born is divested of all quality, there remains the \textit{empty form of the common substratum of all two-ities}. It is this common substratum, this empty form, which is the \textit{basic intuition of mathematics}.''.} In \cite[p.~61-62]{Hersh}, Hersh gives a whole list of different and unclear meanings of this term, but also emphasizes its central role in the philosophy of mathematics: ``1. All the standard philosophical viewpoints rely on some notion of intuition., \ 2. None of them explain the nature of the intuition that they postulate., \ 3. Consideration of intuition as actually experienced leads to a notion that is difficult and complex, but not inexplicable., \ 4. A realistic analysis of mathematical intuition should be a central goal of the philosophy of mathematics.''. However, intuition about our internal world of activity is completely clear, it is intuition in the original sense of that word -- ``immediate awareness'' -- devoid of any ambiguities and misunderstandings. In my opinion, mathematical intuition is precisely this intuition.\footnote{This understanding of mathematical intuition is fundamentally different from Kant's or G\"{o}del's \cite{God}  understanding of intuition.}

Major mathematical models, like the models described above, arise from intuition about our internal activities and organization. It is from these concrete activities that the idea of an idealized world emerges, the world that expands and supplements the internal world of activities.\footnote{The very creation of ideas is largely an unconscious process: ``Most thought is unconscious \ldots inaccessible to direct conscious introspection. We cannot look directly at our conceptual systems and at our low-level thought processes.'' \cite[p.~5]{Lakoff}.} In his book \cite{Lane}, Saunders Mac Lane describes this process on a multitude of examples. The table on page 35 of the book shows a whole list of examples of activities from which certain ideas are born, and from these ideas mathematical concepts and models arise. For example, movements contribute to the idea of change, whose formulation contributes to the concepts of rigid motion, transformation group and rate of change; estimating contributes to the ideas of approximation and closeness, and their formulation contributes to the concepts of continuity, limit and topological space. However, what is the nature of idealized mathematical models? What is the nature of the irrational numbers that Dedekind ``creates'' to fill the ''gaps'' in the linearly ordered field of rational numbers \cite{Ded}? Hilbert writes about infinity as a paradigmatic example of an ideal mathematical element: ``\ldots nowhere is the infinite realized; it is neither present in nature nor admissible as a foundation in our rational thinking \ldots The role that remains to the infinite is, rather, merely that of an idea -- if, in accordance with Kant's words, we understand by an idea a concept of reason that transcends all experience and through which the concrete is completed so as to form a totality \ldots'' \cite[p.~190]{Hi}. In \cite[p.~4]{Fef} Feferman writes:  ``The basic conceptions of mathematics are of certain kinds of relatively simple ideal world pictures that are not of objects in isolation but of structures'' and ``The basic objects of mathematical thought exist only as mental conceptions''. For Hersh, mathematical objects are ``mental objects with reproducible properties'' \cite[p.~66]{Hersh}. 

I consider all the above descriptions of mathematical objects and concepts to be insufficiently clear because they refer to insufficiently clear psychological terms. The same applies to other descriptions of mathematics as a human invention, which I found in the literature. My view of the nature of mathematical models stems from my view of the essential role of language in thinking. My view of the essential role of language in thinking has its inspiration primarily in the works of von Humboldt \cite{Humboldt} and Whorf \cite{Whorf}, and is explained in detail in \cite{Cu0}. According to this view, language is not only a medium for expressing and communicating thoughts, but a medium in which thoughts are realized, a medium in which thoughts take their completed form. In the words of von Humboldt: ``Language is the formative organ of thought. Intellectual activity, entirely mental, entirely internal, and to some extent passing without trace, becomes, through sound, externalized in speech and perceptible to the senses. Thought and language are therefore one and inseparable from each other. But the former is also intrinsically bound to the necessity of entering into a union with the verbal sound; thought cannot otherwise achieve clarity, nor the idea become a concept. The inseparable bonding of thought, vocal apparatus and hearing to language is unalterably rooted in the original constitution of human nature, which cannot be further explained \ldots without this transformation, occurring constantly with the help of language even in silence, into an objectivity that returns to the subject, the act of concept formation, and with it all true thinking, is impossible.'' \cite[p.~50]{Humboldt}. 

It is also important to point out here that,  concerning thinking, the abstractions\label{abstraction} are the language abstractions, and not thinking about the so-called abstract objects. We talk about concrete objects (which can be real or imagined), and abstraction is achieved by extracting certain predicate and function symbols with which we talk about objects. For example, we count using concrete objects. Thus, the language of arithmetic talks about concrete objects (whose nature is not important to us), and with the language itself we achieve the appropriate abstraction. The language of arithmetic separates what is important to us when we use objects as numbers (first number, successor, predecessor, comparison, \ldots) from what is not important (e.g., size of marbles if we use collections of marbles for numbers, or font of decimal digits if we use sequences of decimal digits for numbers). 

In this way of looking at the nature of thinking, which in its final effect is the creation and use of language, we can only realize mathematical models by building an appropriate mathematical language. By choosing names, function symbols and predicate symbols, we shape the initial intuition into a structured conception. Since the conception usually goes beyond our constructive capabilities, the constructed language has only partial interpretation, and that interpretation is internal -- the interpretation in our internal world of activities. Since the interpretation is only partial, and because the imagined domain of interpretation is usually infinite, we cannot determine the truth values of all sentences of the language. Therefore, we must further specify the conception by using an appropriate set of axioms. Thus, the final mathematical model (theory) is a junction of axioms and partial internal interpretation of an adequate language. Sometimes, as we have seen on the example of the group theory, a mathematical model can be reduced to a set of axioms without an internal interpretation, although it arose from a corresponding intuition about the world of internal activities. Sometimes, the internal interpretation can be a total interpretation in another mathematical theory, as we have seen on the example of Euclidean geometry. 

This way of looking at mathematical models can be acceptable to those who believe that mathematical concepts exist even without language. Given the concrete nature of language, unlike the vague nature of mental phenomena, they can accept mathematical language as a convenient representative of mental conceptions. I would also like to point out here that viewing  abstract thinking as the creation and use of language does not mean denying the complex thought processes that take place behind it.  This is a functional view: only the final effect of thinking is considered. 

Furthermore, we must not forget that although a mathematical model is the final product of modelling an intuition about our internal world of activities, in real mathematical practice it is never isolated from the source from which it originated. This is especially important because a mathematical model, generally speaking, is not complete -- there are multiple interpretations that are extensions of partial interpretation and that satisfy axioms; and intuition always leaves room for completion. In addition to testing a mathematical model as a means of rational cognition, it is also tested by how well it models the initial intuition, i.e. how well its consequences correspond to the intuition from which it emerged. In \cite{Lak}, Lakatos clearly demonstrated the importance for mathematics of this internal testing and revision of mathematical models (Lakatos calls this activity quasi-empirical mathematics). Feferman \cite{Fef} distinguishes mathematical conceptions by how close they are to everyday practice. The more distant they are, the less clear he considers them. For those conceptions that are completely clear, such as the conception of natural numbers, he believes that every statement within such a conception has a definite meaning, and thus has a truth value, regardless of whether we are able to determine its truth value or not, so the corresponding logic is the classical first-order logic. For those conceptions at the other end of clarity, such as the set theory, which lack some aspect of definiteness, he believes that the concept of truth may be partial and that the appropriate logic is semi-intuitionistic. Similarly, Ferreirós \cite{FerP} distinguishes mathematical conceptions according to whether their truth is based in our cognition and direct practice, such as the conception of natural numbers (Ferreirós then speaks of elementary mathematics characterized by certainty) or require additional hypotheses, such as the conception of the continuum (Ferreirós then talks about advanced mathematics, which is characterized by the presence of hypothetical statements -- statements that gain or lose their objectivity through mathematical practice). The newer philosophy of mathematics, the so-called philosophy of mathematical practice, devotes special importance to these processes of interaction of mathematical models and basic intuitions, as one of the main sources of mathematics (see, for example, \cite{Man,FerP}).

Although I consider it important how strongly mathematical models are connected to the internal world of activities, my view of mathematical models is more uniform. In my opinion, the basic criterion for evaluating mathematical models is their success as a tool of rational cognition. Seen in this way, there is no difference, for example, between the arithmetic and Hilbert spaces, although their connection to the internal world of activities is different. Also, in my opinion, there is no difference in the meaning of the truth of different mathematical models. To me, all mathematical statements are specifications, whether they specify what we will do with natural numbers in our internal world of activities, or whether they specify the truth of the continuum hypothesis. From such an understanding of the truth of mathematical concepts, it follows that the logic of mathematical models is classical logic, although perhaps in some situations, with regard to the problem that is being solved, it makes sense to look at mathematical models with non-classical logics. 

Just as mathematical models are not isolated from the sources from which they arose, they are not isolated from each other either. Set theory is a natural environment for formulating and comparing mathematical models. In such an approach, axioms become the definition of a certain type of structure. However, set theory analyses the described structures in a uniform way, without going into their nature, whether they are extracted from the external world or from the internal world of our  activities. Thus, although it gives an elegant mathematical description, set theory can also hide the true nature of mathematical models.

Mathematics is largely an elaborated language. The ``magic'' of mathematics is largely the ``magic'' of language. Inferring logical consequences from  axioms, we establish what is true in a mathematical model. This can be very creative and exciting work and it seems that we discover truths about an existing exotic world, but we only unfold the specification. The key difference from scientific theories is that the interpretation here is in our internal world of activities and not in the external world. The external interpretation of a scientific theory enables us to test the theory, whether it has a power of rational cognition of nature. If the theory has such power then at least some objects of the theory  exist in the primary sense of the word and at least some sentences of the theory  are truths about nature. If the language does not have such a part, and that is the case with mathematics, then the objects we are talking about exist only within the conception (story), although they do not exist in the external world. Equally, if the language does not have an experimentally verifiable part then sentences we consider true within the conception are not  true in the external world. We cannot experimentally verify that $||+||=||||$ ($2+2=4$), not because it is an eternal truth of numbers, but because it is the way we add tallies. Likewise, we cannot experimentally verify that $(x^2)'=2x$ because it is the consequence of how we imagined real numbers and functions. Mathematical objects are, possibly, through a partial internal interpretation, objects extracted from our internal world of activities, and mathematical truths are, possibly, through a partial internal interpretation, truths about our internal activities.\footnote{Some statements get their truth through the partial interpretation of language, for example that 2+2 = 4, while for some statements the partial interpretation is not enough to determine their truth  (e.g. the axiom of the completeness of real numbers or the continuum hypothesis). Then the assignment of truth to these statements is, at least in part (the axiom of the completeness of real numbers) or completely (the continuum hypothesis), a specification of an idea that transcends the internal world of activity.}  We are free to imagine any mathematical world. The real external existence  of such a world is not important at all; all that matters is to be  a successful thought tool  in the process of rational cognition. In Cantor's words, ``the essence of mathematics lies precisely in its freedom'' \cite[p.~66]{C}. The only constraint is, inside  classical logic, that conceptions must not be contradictory. For Hilbert, in mathematics to exist  means to be free of contradictions. In Hilbert's words: ``the proof of the consistency of the axioms  is at the same time the proof of the mathematical existence'', \cite[p.~265]{H}. In Dedekind's words, ``numbers are free creations of the human mind'' \cite[p.~vii]{Dedekind}.  

These views are in  sharp contrast with historical views that mathematical truths exist really in some way and that we discover them and not create them. Historically, this change of view  occurred in the 19th century  with the  appearance of non-Euclidean geometries. The new philosophical view of mathematics  has freed the   human mathematical powers and it has caused the  blossoming of modern mathematics. It is a nice example of how philosophical views can influence science in a positive way. According to the old views mathematical truths are a particular kind of  truths about the world. An exemplar is the Euclidean geometry -- according to the old views, it discovers the truths about space. The appearance of non-Euclidean geometries, which are incompatible with Euclidean geometries but are equally logical in thinking  and equally good candidates for the ``true''  geometry  of the world,  has definitely separated   mathematics from  the truths about nature. It has become clear that mathematics does not discover the truths about the world. If it discovers the truths at all they are at best the truths about our own activities in that world. From my personal teaching experience, I know that looking at mathematics as a free and creative human activity is a far better basis for  learning mathematics than looking at it as an eternal truth about an elusive world. Claims that mathematical objects exist in the external world, real or special, and that mathematical truths are truths about such a world, are unfounded and lead to religion. Such a belief can of course be very inspiring and can produce very powerful mathematics, but mathematics itself is a human creation, in its final form a junction of axioms and a partial internal interpretation of an adequate language.

\section{On the existence of mathematical objects}
I will now analyse in more detail the question of the existence of mathematical objects. The key to the answer is in understanding the essential role of language in our thinking. In \cite{Cu0} it is shown how we synthesize our rational cognition of the world through language. An ideal language is, for example, an interpreted language of the first-order logic in which we know the semantic values of all non-logical primitive symbols of the language, as well as the semantic values of all descriptions and sentences in the language. However, in a real process of rational cognition, we use names for which we do not know completely what they name, predicate and function symbols for which we do not know completely what they symbolise, and quantified sentences for which we do not know if they are true or not.  For a theory to be a scientific one, at least some names and some function and predicate symbols must have an exterior interpretation, an interpretation in the exterior world, not necessarily a complete one. This partial external interpretation enables us to perform at least part of the binary experiments described by atomic sentences. This allows nature to put its answers into our framework, so that we can test our conceptions experimentally. Without this part the theory is unusable. On the other hand, due to the partial external interpretation of the language and the impossibility to perform all binary experiments determined by atomic sentences, we necessarily complete the theory with a set of sentences that we consider true in a given situation (axioms of the theory). Thus, scientific theory is also a junction of axioms and partial external interpretation of language. Viewed at the level of the final product, scientific and mathematical theory differ in that the former has an external partial interpretation and the latter an internal one. 

As I have already mentioned when I considered how Euler's number $ e $ exists, although our language usually has only a partial interpretation, the logic of using the language assumes that it is a semantically complete language, i.e.,  that it is an ideal language in the sense as I described at the beginning of this paragraph. Because of this assumption, in thinking itself there is no difference whether we think of objects that really exist or we think of objects that do not really exist. That difference can be registered only in a ``meeting'' with reality. Thus, the question of the existence of the objects we speak and think about is completely irrelevant as long as we do not try to connect language with the external world. Because of the “encounter” with reality, at least some of the objects that scientific theory speaks of must exist in the primary sense of the word, as objects from the external world. 

What about the other objects that scientific theory talks about? Since the early 19th century, physicists and chemists have used the assumption of the existence of atoms to explain many phenomena in matter. Atoms were initially only imagined objects, and in the end, it was established that they really exist. Unlike an atom, ether was initially an imagined object, and in the end, it was established that it does not exist. From this example we see that in science imagined objects can be potentially existing objects. For the next example, let us take a material particle in classical mechanics, the particle occupying a single point in space. It is imagined from the beginning as an idealized object, which does not really exist, but real objects can approximate it. Let us not forget that the basic laws of classical mechanics speak precisely of material particles. E.g. Newton's law of gravitation speaks of the gravitational force between two material particles. Only by applying these idealized laws do we obtain laws of the behaviour of real objects: we ``deconstruct'' real object  into material particles, apply idealized laws to material particles, in order to obtain laws about real objects.  Perhaps the notion of a material point can be avoided in the formulation of classical mechanics, but that formulation, if it exists at all, would be unnecessarily complicated. The above examples show that imagined objects together with statements that we consider true for them are very important for scientific theory. If we were to ban their use, we would literally cripple scientific theories. They are a necessary linguistic tool of scientific theories, and their status can change over time. However, I would point out that their significance follows exclusively from their connection with  the objects of the theory that exist in the external world. 

The situation is analogous to mathematical theories, only here we have a partial interpretation in our internal world of activities. Just as physical models are idealizations and approximations of real processes in the external world, mathematical models are idealizations and approximations of processes in our internal human world of activities. However, as with physical models, the importance of mathematical models also stems from the success of their application in rational cognition. Here the ontological situation is even better than in the case of scientific theories because mathematical theories are completely under our control, in the sense that we determine what will be true and what will not be true in such a theory, not the external world. Unlike objects from the outside world, which are under the authority of nature, we create mathematical objects, and the truths about them are in fact their specifications. Some mathematical objects are determined up to isomorphism, such as natural numbers. We realize these objects internally as part of the (partial) realization of the corresponding structure. Some mathematical objects can be completely realized, such as not too large natural numbers, or up to a satisfactory approximation, as irrational numbers. Some mathematical objects have a more specific nature, such as sets. We realize them through their close representatives, for example, sets with appropriate one-place predicates (linguistic forms of a certain language). Representation is a kind of isomorphism. However, the difference is that not all isomorphic structures are equal to us, but we have a certain, perhaps not entirely clear interpretation as with sets, to which we look for ``close'' representatives. Some mathematical objects have a completely determined nature, such as geometric objects in the primary interpretation of Euclidean geometry. We realize them approximately, but directly (not through isomorphism or representation) in the world of our internal activities.

There is another way we can understand the existence of mathematical objects. If a theory is consistent, from Henkin's \cite{Hen} proof of the completeness of a deductive system follows the existence of a canonical model whose objects are classes of equivalence of the corresponding terms of the language. Since this model is homomorphically embedded in any interpretation that satisfies axioms, this means that we can always represent a finite portion of imagined mathematical objects through their proper names in the language in which we speak of those objects.

To summarize, let us start from the fact that the question of the existence of the objects of language is irrelevant to the formation and use of language in the process of thinking, until the moment when we apply language to the external world. For a scientific theory, at a given stage of its development, the existence of parts of language that speak of objects whose existence in the outside world has not been established and about which the theory makes certain claims can only be evaluated through the way these parts relate to the experimentally verifiable part of the theory. Through this interaction with the experimentally verifiable part of the theory, it can be shown that such objects either exist or do not exist, and thus their ontological status is revised. However, for some objects, such as material particles, it can be shown that they are imagined objects that enable an efficient, perhaps necessary, linguistic synthesis of rational cognition. In the latter case, we may regard them as mathematical objects, our imagined tools for rational cognition. Thus, we can equate them with the parts of language that belong to mathematics and reduce the problem of their existence to the problem of the existence of mathematical objects. Since mathematical language, i.e., mathematical models, do not have an experimentally verifiable part, from that point of view, the question of the existence of mathematical objects is irrelevant, and mathematical truths are only the specifications of imagined objects.\footnote{The specifications are rarely complete in the sense that they determine the truth value of each statement.}  The only thing that matters is how we apply mathematical language in rational cognition. Thus, in the broadest view, a view that does not restrict mathematical freedom, a mathematical model is a set of sentences of a language that we declare to be true and that has possibly a partial internal interpretation. Thereby, it is important that the mathematical model is intended to be a successful  tool for rational cognition or rational activities in general. It is the success of a mathematical model that determines its quality. And for a mathematical model to be successful, its language must somehow be related to the language we use to describe the outside world. For example, when we count, basically, with the words `` one '', `` two '', `` three '', \ldots, we associate the objects we are counting. We do not need to know at all what those words name and whether they name something at all. Of course, we imagine the language of numbers so that these words name objects in some structure of numbers (we can assume that these words themselves form that structure). But the point is that we can count without knowing what those words name. Let us imagine another situation, that someone gave us axioms of real numbers without us ever hearing about real numbers. By deriving statements from these axioms and defining new words using existing ones, we may understand that this language allows us to describe the measurement process and express the measurement results in it. In this way, that language (that theory) itself  becomes a successful tool of rational cognition, regardless of what the objects of the language are and whether they exist at all. I have given these extreme cases to show that in the broadest view of a mathematical model, as a consistent set of axioms in a language, the question of what that theory is talking about does not need to be important at all for its use. The set of axioms can even have multiple interpretations.  And if you just want to have an interpretation, we can build a canonical model from the very strings of the language of the theory. But real and successful mathematical models do not arise in such a way. As I described in this article, they are the results of modelling intuition about activities from our internal human world -- that world gives life and power to a mathematical model and in that world the mathematical  model has its own partial interpretation which significantly reduces the incompleteness of the model. To conclude, the best I can say is that mathematical objects do not exist in the external world -- they are our internally imagined objects, some of which, at least approximately, exist in our internal world of activities or we can realize or represent them there.

\section{How mathematics is applied}
Here I will explain in detail how is it possible that something imagined can contribute to  rational cognition of the world? At the beginning of Section 7 I stated that we can use mathematical models  directly or indirectly. In the same paragraph I also stated how mathematical models are used indirectly. I think their indirect use is clear enough, so in this section I will only deal with the direct use of mathematical models.

My understanding of the direct use of mathematical models stems from my view of rational cognition as a kind of synthesis of us and nature that takes place through the use of language.  Let us look at the simplest statement about the world of the form P(a) where ``a'' is the name of an object and ``P'' is a predicate symbol. For example, ``a'' is my dog's name and ``P'' stands for ``is afraid of thunder''. The structure of the statement P(a) reflects our innate approach to the world which we divide  into objects (elements  upon which  something is done) and into predicates (which  determine what is done). To my knowledge, Whorf is the first one to recognise that the object-predicate dualism is a prominent feature of Indo-European languages: ``Our language thus gives us a bipolar division of nature. But nature herself is not thus polarized.'' \cite[~p.247]{Whorf1}. He also recognizes that the dualism and  the way we analyse nature is not inherent to nature but to our approach to nature: ``We dissect nature along lines laid down by our native language. The categories and types that we isolate from the world of phenomena we do not find there because they stare every observer in the face; on the contrary, the world is presented in a kaleidoscopic flux of impressions which has to be organized by our minds -- and this means largely by linguistic systems in our minds. We cut nature up, organize it into concepts and ascribe it significance as we do, \ldots '' \cite[~p.231]{Whorf1}. To determine whether a statement P(a) is true, for example whether my dog is afraid of thunder, knowing the meaning of the symbols ``a'' and ``P'' is necessary but not sufficient. I still have to do an appropriate experiment, let nature give its contribution, to determine that it is a true sentence. These binary investigations (``experiments'') are the starting point for the overall rational cognition. We pose a question and offer two possible answers, the so-called truth values termed  $True$ and $False$,  and  nature selects an answer. The selected truth value does not belong exclusively to us, nor does it belong exclusively to nature. It is the objective result of the synthesis of us and nature in the process of rational cognition: it  discriminates what is and what is not. This synthesis is also extended to more complex linguistic forms, according to the meaning we give to these forms. The truths that we achieve   are not  truths about the world itself -- they are truths of our rational interaction with the world. In this process we form a rationalized reality. The rationalized reality is the result of our synthesis with nature through the creation and use of language in the process of rational cognition. As Whorf writes: \cite[~p.285]{Whorf3}: ``We don't think of the designing of a radio station or a power plant as a linguistic process, but it is one nonetheless.''  A detailed presentation of my view on the role of language in rational cognition can be found in \cite{Cu0}.

The role of mathematics in the above-described synthesis of us and nature in the process of rational cognition is precisely the one I stated in the second paragraph in Section 7: mathematics is an inner organization of rational cognition and knowledge, a thoughtful shaping of the part of the cognition that belongs to us. Using the example of real numbers, natural numbers, and a first order language, I will show how mathematical models are directly used in rational cognition understood in this way.
 
Real numbers are imagined as the results of the process of measuring. In the process of measuring, we connect them with the external world, enabling nature to select one of the offered numbers as its answer. The number itself is not real (in the sense that it does not belong to the external world) but nature's selection is real. Numbers belong to our inner organization of the measurement process but nature's selection is a truth about the external world.  For example, during the process of measuring the speed of light, among all other numbers, nature  selects the number $c$. The selected number $c$   possibly exists as our internal construction. Whether it is a rational or irrational number depends on the choice of units of measurement. However, that $c$ is the speed of light is the idealized truth about the external world which is synthesized in the process of measuring. Idealized, because we assume that c is the result of an idealized process of measurement to which the actual measurement is only an approximation. 

 The simple assertion about natural numbers, that $2+2 = 4$, is a true sentence about the imagined world of natural numbers, and not a truth about the external world. However, through the real process of counting,  we can use  assertions about numbers to obtain  synthesizing  assertions about the external world. For example, when we put two apples in a basket which already contains two apples,  we predict that there will be   $2+2=4$ apples in the basket. This is the prediction about reality deduced from the mathematical assertion that $2+2=4$ and the assumption that  the mathematical model of counting and adding one set of things to another set  is applicable in this situation. However, we must distinguish the mathematical assertion that $2+2=4$ from  the assertion about  reality that when we add $2$ apples to $2$ apples there will be $4$ apples. The best way to see the difference is to imagine a situation where we add 2 apples to 2 apples  and get 5 apples. It would mean that it is not always true, as we have thought, that adding 2 apples to 2 apples gives 4 apples, but in some situations, according to  as yet unknown physical laws,   an additional apple emerges. However, this situation would not have any influence on the world of numbers. In that world it is still true that $2+2=4$. It only means that in some real situations we cannot apply the mathematical model of counting and addition. 
 
 The natural numbers model of counting, as well as the real numbers model of measuring or any other mathematical model, have their assumptions of applicability. For the natural numbers model, we assume that we can associate the number of elements to a collection of objects by the process of counting. For the real number model,  we assume that we can always continue to measure with a ten times smaller unit, if it is necessary. A real process of measuring, for example, the distance of point $B$ to point $A$, must stop in a certain step, because the passage to a ten times smaller unit would not be possible with an existing measuring instrument or that passage revises our understanding of what we are  measuring at all. For example, if we measure the distance between points A and B marking it with a pencil on a paper, the passage to the one hundredth part of millimetre requires a microscope. Looked at under a microscope, $A$ and $B$ are not  points any more, but  diffused flecks. And what are we measuring now? The distance between the closest points of the two flecks? If we continue to magnify, we will see molecules which constitute flecks and which are in constant  vibrations. And if we moved on to even tinier parts, we would come to the world of quantum mechanics in which  classical  notions, on which our conception of measuring  distance is based, are not valid any more. However, the question whether it is possible to apply the mathematical model of measuring in reality  is not a mathematical question at all.       Only when the assumptions of the model of  measuring are (at least approximately) fulfilled can we apply the mathematics of real numbers to the real world. Likewise, only when the assumptions of the model of counting are fulfilled can we employ the mathematics of natural numbers  to the real world. Only then we have at our disposal the whole mathematical world that can help us in asserting truths about the real world.  We have at our disposal an elaborate non-verifiable language which we can connect with a verifiable language, mathematical truths which we can synthesize  into  truths about the external world. If contradictions occur in an interaction with the verifiable part of the language, it does not mean that the mathematical model is false (the concept of the real truth and falsehood does not make any sense for the model), but that the assumptions about its applicability in that situation are false. 

A first-order language is a mathematical model constructed for the use in rational cognition just like natural numbers are constructed for counting.\footnote{See, for example, \cite{vanD} for a detailed description of the syntax and semantics of first-order languages.} This model is the result of   thoughtful modelling of intuition about our natural language. A first-order language has the external assumptions of its use. These are: 1) existence of the domain of interpretation, 2) every name names an object, 3) every function symbol symbolizes a function which applied to objects gives an object, and 4) every predicate symbol symbolizes a predicate which applied to  objects gives one of the two possible results, ``true'' or ``false''. Only when these assumptions are fulfilled can we employ the first-order language  model  to the real world: all the first order logic together with computability theory to examine external structure and derive truths about it.\footnote{In my reading of Fenstad \cite{Fen}, this is exactly the general scheme of meaning and application of mathematics that he proposes.}

Quine \cite{QuineN} in his naturalized epistemology considers that every part of the web of knowledge is liable to experiment, including logic and mathematics. That is true, but there are qualitative differences between science on one side and logic and mathematics on the other. Experimental evidence can affect the truth values of scientific sentences but not the truth values of mathematical and logical sentences. It can only   question the applicability or adequacy of mathematical models and language frameworks in some parts of science. Scientific theories are true or false of something while mathematical models are good or bad of something.

%The nature of our thought and use of language leads to extracting a certain structure from the subject of thinking. We extract from the subject certain objects, relations and operations and we describe their properties. From the language point of view, we name certain objects, we introduce certain predicate and function symbols and consider certain sentences which are built from these symbols to be true. This process of abstracting enables us to successfully manage  a vast complexity of the world on one side, and it binds us to the extracted structural framework on the other side. In that  way we cannot cognise objects ''in themselves'' but only how  they are connected mutually. Our rational cognition of the world is necessarily structural. However, as opposed to external structures, the structures  that are supported by reality,  where abstraction is the result of abstracting  elements from reality, mathematical structures are imagined structures  that ''emerge'' from our inner world and abstractions associated with them   are mainly without support. They are the descriptions without the described. If a mathematical abstraction has  support, the support is of our human nature, it consists of symbols, actions, etc.  From the language point of view, mathematical theories are mainly sets of sentences without interpreted parts. If they have an interpreted part, it is  an inner interpretation, in  another mathematical theory or directly in our human world. 

\section{The structure of mathematics}
Various mathematical models are not mutually disconnected but they are interwoven. Moreover,  we  express these connections also by corresponding mathematical models. First of all, there is a not so big collection of primitive mathematical models (``mother structures'' in Bourbaki's terminology \cite{Bourbaki}) that model the basic intuitions about our internal world of activities, intuition about near and remote (topological and metric structures), about measuring (spaces with measure), about straight and flat (linear spaces), about symmetry (groups), about order (ordered structures) etc. We use them  as ingredients of more complex mathematical models. The complex mathematical models   enable us to realize some simple and important mathematical ideas (for example, we use normed linear spaces to realize an idea of the velocity of change) or they have  important applications (like Hilbert spaces which,  among other things, describe the states of quantum systems). Therefore,  the world of mathematics is built   of some primitive models which model our ideas about our internal world of activities and of various ways of comparing and combining  these models into more complex models. Corresponding mathematical theories are interpreted mutually or have  internal interpretations. If a theory has an external interpretation then only its sentence part belongs to mathematics. For example, classical mechanics is a theory about the world. However, if we ignore its externally interpreted part, we get a set of sentences which we can investigate mathematically. Hence, we can say that mathematics is concerned with the internal models and internal properties of external models, or more simply, it is concerned with that part of  rational cognition that belongs to us. Furthermore, such diverse mathematical models are the basis for secondary mathematical models that model how to compare  structures (set theory and category theory) and in what language to describe them (mathematical logic). 

However, regardless of the complexity of the world of modern mathematics, mathematics is an inner organization of rational  cognition based on the modelling of the intuition about our internal world of activities. In constructing mathematical worlds, the criteria of real truth and falsehood have no meaning, although every such  world has its inner truths and falsehoods which shape and express the underlying conception. However, in the process of rational cognition we synthesize mathematical objects and truths into truths about the real world (for example, into Newton's law of universal gravitation). In the construction of mathematical worlds, the criteria of simplicity and beauty, almost artistic criteria, are of real importance as well as how well  the constructed theories model the original intuition and ideas. The fulfilment of these criteria is a good indicator, as experience shows,  that the main criteria will also be fulfilled, the criteria of the direct or indirect usefulness of those worlds as our thinking  tools in the process of rational cognition. If the ideas are good and if they are  well modelled mathematically, sooner or later they will certainly find a successful application, as we have seen in the example of Riemann's manifolds.

\section{Some tests for the view of mathematics as an imagined tool for rational cognition}

I consider that this view of mathematics satisfies both concerns of  Benacerraf dilemma \cite[p.~661]{Be}: ``(1) the concern for having a homogeneous semantical theory in which semantics for the propositions of mathematics parallel the semantics for the rest of the language, and (2) the concern that the account of mathematical truth mesh with a reasonable epistemology''. I have shown in Section 8 that the only difference between mathematical and scientific language is whether it is a partial interpretation of language in the internal world of activities or in the external  world. Since we are the creators of the mathematical worlds, their epistemological status is unquestionable.

Also, I consider that the analysis conducted in this article has shown that this  view of mathematics  fully meets Bueno's  \cite[p.~63]{OB} ``five desiderata that an account of mathematics should meet to make sense of mathematical practice (my brief answers are in parentheses): (1) The view explains the possibility of mathematical knowledge (yes, it does -- we create mathematical knowledge). (2) It explains how reference to mathematical entities is achieved (yes - they are imagined objects, some of which, at least approximately, exist in our internal world of activities or we can realize or represent them there). (3) It accommodates the application of mathematics to science (yes -- this is most explicitly shown in Section 9). (4) It provides a uniform semantics for mathematics and science (yes -- this is explained in the paragraph above). (5) It takes mathematical discourse literally (yes -- this is most explicitly explained in Section 9 when describing what a mathematical model is).''

Concerning Quine-Putnam indispensability argument for the existence of mathematical objects as it is spelled in \cite{Co},  this view of mathematics  fully supports the second premise stating that ``mathematical entities are indispensable to our best scientific theories'' and rejects its first premise, which claims that ``we ought to have ontological commitment to all and only the entities that are indispensable to our best scientific theories''.  I have  shown above in the simplest example of counting, which is certainly indispensable to our best scientific theories, that one can use the language of numbers to count  in (scientific) application without even knowing what numbers are and whether they exist at all.

\section{Mathematics as an imagined tool for rational cognition and other views of mathematics}

Since I am a mathematician by profession, it was difficult for me to find suitable philosophical terms to explain my view of mathematics. It is even more difficult for me to compare my view with the multitude of existing philosophical views on mathematics. As far as I am aware, the view of mathematics presented in this paper has a certain originality, if nothing else, because on the one hand it elaborates some insufficiently elaborated views of mathematics, and on the other hand, it encompasses some prominent views of mathematics but, as a whole, avoids their one-sidedness.

It is clear that this view of mathematics has nothing in common with realistic views of mathematics, according to which mathematical objects and mathematical worlds belong to the external world. In my opinion, such views are full of metaphysical ghosts to which I cannot add any value, and which, on the other hand, can be dangerous to rational cognition, in the same way that religion can be dangerous.\footnote{I accept that such views can be inspiring and motivating for one's mathematical activity, but I cannot accept them as truths about mathematics.} In the view of mathematics described in this paper, the human being and the human community create mathematics, just as they create, for example, works of art. 

This view of mathematics is close to Hersh's humanistic philosophy of mathematics and Ernest's social constructivist philosophy of mathematics and can be considered a certain elaboration of their views in that part where their views are insufficiently developed. According to them, the basis of all mathematics is real mathematical activity, individually performed, and socially supported and maintained.  In \cite[p.~63]{Hersh} Hersh writes ``A world of ideas exists, created by human beings, existing in their shared consciousness. These ideas have objective properties, in the same sense that material objects have objective properties.'', and that  ``any mathematical object you like — exists at the social-cultural-historic level, in the shared consciousness of people (including retrievable stored consciousness in writing). In an oversimplified formulation, "mathematical objects are a kind of shared thought or idea."''. In \cite{Ernest} Ernest writes: ``Social constructivism adopts an approach to mathematical objects perhaps best described as nominalist, regarding them as linguistic/conceptual objects.'' (page 261), and in more detail, ``According to the social constructivist view the discourse of mathematics creates a cultural domain within which the objects of mathematics are constituted by mathematical signs in use. Mathematical signifiers and signifieds are mutually interacting and constituting, so the discourse of mathematics which seems to name objects outside of itself is in fact the agent of their creation, maintenance, and elaboration, through its use.'' (page 193). From the above, it is clear how important the linguistic component of mathematics is in Ernest's approach, as well as in mine. However, I consider Hersh's and Ernest's descriptions of mathematical objects to be insufficiently clear. Perhaps the notion of a mathematical model developed in this paper could serve as a clarification of their views on mathematical objects and concepts. On the other side, the humanistic view  of mathematics complements the view of mathematics described in this paper. 

Likewise, the embodied mind and embodied cognition approach to mathematics, as they are vividly described in \cite{Lakoff}, complement the view of mathematics described in this paper. These approaches show the importance of the biological component of mathematics and show the connection of the mathematical ideas with our biological predispositions and activities. Research in this area seriously dis-confirms realistic conceptions of mathematics and confirms the thesis that mathematical ideas are deeply rooted in our behaviour and that mathematics is created rather than discovered. Likewise, Kitcher's historical explanation of mathematics \cite{Kit} and emphasis on the importance of mathematical practice \cite{Man,FerP} complement my views. These are very important aspects of mathematics that shape current mathematical conceptions but I do not find in these works an elaborate view of what mathematical conceptions are. All those who think, like me, that mathematics is a human invention, face the same problem -- what kind of human invention is mathematics? The notion of a mathematical model explained in this paper is the answer I came up with.

I consider that the view of mathematics presented in this paper belongs in large part to the tradition which has begun with Dedekind, Cantor and Hilbert. The creationist and structuralist views of Dedekind about mathematics have spread  into mathematics community through the works of Emmy Noether, van der Waerden  and Bourbaki, and they have become the trademark of modern mathematics. I was surprised when I saw the amount of attention paid by the philosophical community to the Benacceraf’s article titled ``What numbers could not be'' (\cite{Benacerraf}), because I believed that the contents of the article were well known from the time of Dedekind. The subject of discussion is whether and in what form psychologism is present in Dedekind's work, also in the quotation from the beginning of the article (see \cite {Re} for a detailed analysis). If Dedekind, under the free creation of natural numbers (as well as other numbers), considered their creation as a creation of abstract thoughts then that is different from  the view of mathematics set out here. In this view, mathematical objects are always imagined as concrete objects, because that is how our language works, regardless of whether we can realize them in our internal world of activities or not. We already do abstraction with language itself, which ensures that we consider only those properties of mathematical objects that interest us. Thus, language is the bearer of the required abstraction  and not the objects of language.

In the view of mathematics presented in this paper, the final mathematical model (theory) is a junction of axioms and partial internal interpretation of an adequate language. Thus, this view encompasses formalism (in the part related to axioms) and constructivism (in the part related to partial internal interpretation). The subject of formalism is formal systems, and the subject of constructivism is the internal world of our activities. But if we do not look at them as methodological approaches but as philosophies, they are one-sided -- formalism limits mathematics to formal manipulations of language forms and constructivism to the internal world of activities.

Logicism, regardless of its success or failure, or its modern reincarnation in classical set theory, is also a one-sided view of mathematics. In my view of mathematics, logicism and set theory are only large (and important) mathematical models in which  `` all '' mathematical models can possibly be represented up to an isomorphism but which neglect the existence of particular interpretations. For example, in set theory we can use ordered triples of real numbers to realize the structure of three-dimensional Euclidean space but we cannot reconstruct Euclidean geometry in its primary interpretation, as the structure of the space of our internal activities.
 
This view of mathematics can fit into naturalism, if we understand it in the broader sense of the word,\footnote{In \cite{nat} one can see how many different variants of naturalism there are.}  because human cognition is basically holistic and every part, even mathematics, can be properly perceived only in relation to that whole. However, it cannot fit into naturalism in the narrower sense of the word which sets ontological conditions in the first premise of the Quine-Putnam indispensability argument and does not respect the difference between logic and mathematics on the one hand and the natural sciences on the other, as I commented at the end of Section 11.

 The structuralism that is part of this view of mathematics is structuralism that  is present in mathematical practice. It is deprived of all philosophical additions on the nature of mathematical structures and mathematical objects by which various structuralist views in the philosophy of mathematics are distinguished, as nicely dissected in \cite {Re1}. In what follows, by structuralism I will mean precisely this structuralism. The approach in mathematics according to which only the structure of objects is important is named in \cite {Re1} structuralist methodology and is described as follows: ``Mathematicians with a structuralist methodology stress the following two principles in connection with them: (i) What we usually do in mathematics (or, in any case, what we should do) is to study the structural features of such entities. In other words, we study them as structures, or insofar as they are structures. (ii) At the same time, it is (or should be) of no real concern in mathematics what the intrinsic nature of these entities is, beyond their structural features.''. I consider that structuralism is only one, although very important, aspect of mathematical modelling but we cannot reduce all mathematics to it. As I stated above, the nature of our thought and use of language, as well as the way in which we navigate the vast complexity of the world, leads to extracting a certain structure from such a domain. Thus, in the study of a phenomenon, we limit ourselves to the study of a certain structure of that phenomenon. However, structuralism is an approach that studies structures formally and does not enter into their nature. It thus naturally falls under mathematics, and the language of set theory is a natural language for the study of structures. However, not only do the structures we extract from the outside world (and thus do not belong to mathematics) have their content, but mathematical structures can also have their content. We have seen that only structural properties are important for natural numbers, and here the whole mathematics is covered by the structural approach. However, what about sets, for example? Although we do not have a completely clear interpretation of what sets are, we do not think that only structure is important for sets, but that they have additional content. The situation is even more pronounced with Euclidean geometry which I consider to be a mathematical model of the space of our human activities, and which is on an equal footing with number systems, and is not a part of physics. Its objects are completely determined in our internal world of spatial activities and are not just ``places'' in a certain structure. Although we can study various structures of the same type as the structure of Euclidean geometry, or connect it with isomorphic structures, thus we cannot exhaust the mathematical content of Euclidean geometry. So, if we do not look at structuralism as a methodological approach but as a philosophy, it is one-sided, or as Machover says in \cite{Mac}, it is ``alienated''.

The main difference between fictionalism and this view of mathematics is that fictionalism considers that mathematical objects do not exist (or at least that it does not matter whether they exist at all) and that mathematical claims are not true, while I consider that mathematical objects are imagined objects that we can at least partially realize in the world of our internal activities and that mathematical truths are specifications of mathematical ideas, possibly idealized truths about the  world of our internal activities. This connection between mathematics (mathematical objects and mathematical truths) and our  world of internal activities is crucial -- it is the source from which mathematics arises and the environment in which it is applied. Furthermore, I believe that fiction is always linguistic fiction. And  I do not think that when we use mathematical language, we pretend that there are mathematical objects and pretend that what we say about them is true. The very nature of language use requires us to assume the existence of the objects we are talking about. Only when we step out from that language can we question the existence of these objects (in this new language), as Carnap explained long ago in \cite{CarnapE}. In my opinion, fictionalism has two different versions. In one version, fictionalism does not take mathematical language literally but figuratively, according to Field \cite{Fi} as a conservative extension of a content language, while Yablo \cite{Ya} considers it as representational aid. The second version (Balaguer \cite{Ba}, Leng \cite{Le}, Bueno \cite{OB} takes the mathematical language literally. Since I, too, take mathematical language literally, my view of mathematics may have common points only with the second version of fictionalism, and I will here refer to this version.  I have  shown above by the example of natural and real numbers how a mathematical model can be applied purely linguistically without even knowing  what its objects are and whether they exist at all. This way of using mathematics is in line with fictionalism and it shows that, if we look at the final product, fictionalism is very close to formalism. However, such an approach to mathematics is artificial and limiting - neither mathematics is practised in this way, nor can all mathematical activities be reduced to such an approach. It is a partial internal interpretation that gives and sustains life to mathematical models.

 .

\bibliography{Mathematics3}
\bibliographystyle{apalike}

\end{document}